\documentclass{research}
\usepackage{mathrsfs}
\usepackage{graphicx}
\RequirePackage{amsmath}
\RequirePackage{amssymb}

\newcommand{\df}[2]{\displaystyle\frac{#1}{#2}}
\newcommand{\tf}[2]{\textstyle\frac{#1}{#2}}

\newcommand{\Intt}[1]{\displaystyle\int_{#1}}

\newcommand{\PDD}[2]{\df{\partial #1}{\partial #2}}
\newcommand{\PD}[1]{\df{\partial }{\partial #1}}

\newcommand{\PDDT}[2]{\df{\partial^{2} #1}{\partial #2^{2}}}

\newcommand{\os}[1]{\overline{#1}}
\newcommand{\oss}[2]{\overline{#1}_{#2}}

\newcommand{\DDt}[1]{\df{D #1}{D t}}

\newcommand{\eps}{\varepsilon}
\newcommand{\B}[1]{\mbox{\boldmath $ #1 $} }

\newcommand{\grad}{\B{\nabla}}
\newcommand{\dive}{\B{\nabla .}}
\newcommand{\curl}{\B{\nabla\times}}

\newcommand{\D}[1]{{#1}^{\prime}}
\newcommand{\DD}[1]{{#1}^{\prime\prime}}

\newcommand{\be}{\begin{eqnarray}}
\newcommand{\en}{\end{eqnarray}}

\newsavebox{\astrutbox}
\sbox{\astrutbox}{\rule[-5pt]{0pt}{20pt}}

\begin{document}

\title[Variational principle for Navier-Stokes equations] {Variational Principle for Velocity-Pressure Formulation of Navier-Stokes Equations}

\author[S.G. Sajjadi]{Shahrdad G. Sajjadi$^\dag$}

\affiliation{Department of Mathematics, Embry-Riddle Aeronautical
University,\\ 600 S. Clyde Morris Boulevard, Daytona Beach, FL
32114-3900, USA.}

\footnotetext{$^\dag$Also visiting scholar at Trinity College, Cambridge.}

%\pubyear{2007}
%\volume{538}
%\pagerange{1-30}
%\setcounter{page}{1}

\label{firstpage}
\maketitle

\begin{center}{\bf Abstract}\end{center}
\begin{abstract}

The work described here shows that the known variational principle for the
Navier-Stokes equations and the adjoint system can be modified to produce a 
set of Euler-Lagrange variational equations which have the same order and same
solution as the Navier-Stokes equations provided the adjoint system has a
unique solution, and provided in the steady state case, that the Reynolds
number remains finite.
\end{abstract}
\section{Introduction}
Given a differential equation subject to some boundary conditions, it is an
important question to decide if a variational principle exists, especially if 
this variational principle is to be used for providing an approximate solution
to the original problem. For some problems minimum and maximum principles may 
hold, leading to reciprocal variational principles which provide a means for
obtaining upper and lower bounds on a variational integral. In other problems,
however, the variational principle, if at all exists, may only be a stationary
principle and no minimum or maximum can be achieved.

For nonlinear and non self-adjoint operators the construction of variational
principle is not well understood. Millikan [1] and
Finlayson [2] have shown that there is no variational principle for the
Navier-Stokes equations involving soley the velocity vector $\B{u}$ and the
pressure $p$ unless either $\B{u\times}\curl\B{u}=0$ or $\B{u.}\grad\B{u}=0$,
when the equations effectively becomes linear.

Under such circumstances we can express the velocity vector
\begin{eqnarray}
\B{u}=\grad\psi+\lambda\grad\mu\nonumber
\end{eqnarray}
in which $\psi$, $\lambda$ and $\mu$ are scalar functions of $\B{x}\in\Omega$
and $t$, $\Omega\subseteq\Re^{3}$. The scalar functions $\lambda$ and $\mu$
satisfy
\begin{eqnarray}
\DDt{\lambda}=\DDt{\mu}=0\label{eq1}
\end{eqnarray}
where $D/Dt\equiv\partial_{t}+\B{u.}\grad$. Except for an arbitrary function of 
$t$ alone, $\psi$ is then determined uniquely by $\lambda$ and $\mu$ {\em qua}
functions of the space variables $x_{j}$ for each $t$. With $\lambda$ and $\mu$
given in the domain $\Omega$, the residual equation of mass conservation
\begin{eqnarray}
\begin{array}{lll}
\dive\B{u}=\Delta\psi+\lambda\Delta\mu+\grad\lambda\B{.}\grad\mu=0 &
{\rm in} & \Omega \\
\B{n.}(\grad\psi+\lambda\grad\mu)=0 & {\rm on} & \partial\Omega \\
\end{array}\nonumber
\end{eqnarray} 
(if $\Omega$ is taken to be bounded rigidly or
\begin{eqnarray}
\begin{array}{lll}
|\grad\psi+\lambda\grad\mu|\rightarrow 0 &{\rm as} & |\B{x}|\rightarrow\infty\\
\end{array}\nonumber
\end{eqnarray}
if $\Omega$ is unbounded), comprise a Neumann problem for $\psi$. The problem  
is soluble in at least a weak form if merely 
$\lambda\grad\mu\in L^{2}(\Omega)$, where $\lambda\in L^{\infty}(\Omega)$ and
$\mu\in H^{1}(\Omega)$. The respective variational principle of $\psi$ is given
by
\begin{eqnarray}
J(\psi,\lambda,\mu)=\min_{f\in H^{1}(\Omega)}J(f,\lambda,\mu)\label{eq2}
\end{eqnarray}
where
\begin{eqnarray}
J(f,\lambda,\mu)=\Intt{\Omega}\left(\tf{1}{2}|\grad f|^{2}+\lambda\grad\mu
\B{.}\grad f\right)d\B{x}+\Intt{\partial\Omega}f\lambda\grad\mu\B{.n}\,dS
\nonumber
\end{eqnarray}

With $\psi$ thus determined as a functional transformation of $\lambda$ and
$\mu$, the Cauchy problem (\ref{eq1}) with solution $(\lambda,\mu)(t)$ appears
complete, and the following conclusion can be drawn. If potentials 
$(\lambda,\mu)(0)$ match the given initial data for the inviscid incompressible
fluid equations, namely the Euler equations, then the evolution of 
$(\lambda,\mu)$ according to (\ref{eq1}) with $\psi$ fixed by (\ref{eq2}) at 
each $t\geq 0$, fully reproduces the solution of the hydrodynamic problem.

In this work we present a variational principle for the Navier-Stokes equations
and the adjoint system. This variational principle is then modified to produce
a set of Euler-Lagrange variational equations which have the same order and
solution as the Navier-Stokes equations provided the adjoint system has a
unique solution.
\section{A variational principle for a non self-adjoint problem}
One of the simplest non self-adjoint boundary value problems is that for a
damped harmonic oscillator with fixed end conditions,
\begin{eqnarray}
\DD{y}+2a\D{y}+by=0,\hspace*{0.5cm}0<x<1,\hspace*{0.5cm}(\,\,\D{)}\equiv 
\df{d}{dx}\label{eq2.1}
\end{eqnarray}
where $a$ and $b$ are constant and $y(0)=\alpha$, $y(1)=\beta$. This problem 
has a unique solution provided $\sqrt{(b-a^{2})}$ is not an integral multiple 
of $\pi$. Finding a variational principle for this problem, involving soley
$y$ and its derivatives without first introducing a transformation of
variables of some type, is not possible. It is easy to show that (\ref{eq2.1})
is exactly equivalent to the variational problem of finding stationary values
of a functional of two independent functions of $x$. The appropriate problem is
the following:

Find stationary values for the functional
\begin{eqnarray}
J(y_{1},y_{2})=\df{1}{2}\int_{0}^{1}\left[(y_{1}^{\prime^{2}}-
2ay_{2}^{\prime}y_{1}-by_{1}^{2})-
(y_{2}^{\prime^{2}}-2ay_{1}^{\prime}y_{2}-by_{2}^{2})\right]\,dx
\label{eq2.2}
\end{eqnarray}
\begin{eqnarray}
b-a^{2}\neq m^{2}\pi^{2},\hspace*{0.5cm}\mbox{$m$ an integer}\label{eq2.2a}
\end{eqnarray} 
among the class of functions
$y_{1}(x)$, $y_{2}(x)$ which have continuous second derivatives and satisfy the
following boundary condtions as $y$,
\begin{eqnarray}
y_{1}(0)=y_{2}(0)=\alpha,\hspace*{1cm}y_{1}(1)=y_{2}(1)=\beta\label{eq2.3}
\end{eqnarray}

The Euler-Lagrange equations for the functional (\ref{eq2.2}) are
\begin{eqnarray}
\left.\begin{array}{ll}
y_{1}: & y_{1}^{\prime\prime}+2ay_{2}^{\prime}+by_{1}=0\\
y_{2}: & y_{2}^{\prime\prime}+2ay_{1}^{\prime}+by_{2}=0\\
\end{array}\right\}\label{eq2.4}
\end{eqnarray}
so that the difference function $\os{y}=(y_{1}-y_{2})/2$ satisfies the
homogeneous problem
\begin{eqnarray}
\os{y}^{\prime\prime}-2a\os{y}^{\prime}+b\os{y}=0,\hspace*{1cm}
\os{y}(0)=\os{y}(1)=0\nonumber
\end{eqnarray}
When the condition (\ref{eq2.2a}) is imposed, which is precisely the condition that
the original problem should have a unique solution, it is easy to see that 
$\os{y}=0$, and that equations (\ref{eq2.4}) with boundary conditions 
(\ref{eq2.3}) become
\begin{eqnarray}
\DD{y}+2a\D{y}+by=0,\hspace*{0.5cm}y(0)=\alpha,\hspace*{0.5cm}y(1)=\beta
\nonumber
\end{eqnarray}
the original problem.

The Lagrangian in the functional (\ref{eq2.2}) may be written in term of
\begin{eqnarray}
y=\df{y_{1}+y_{2}}{2},\hspace*{1cm}\os{y}=\df{y_{1}-y_{2}}{2}\nonumber
\end{eqnarray}
and the functional then has the form
\begin{eqnarray}
J(y_{1},y_{2})&=&2\int_{0}^{1}\left(\D{y}\os{y}^{\prime}+2ay\os{y}^{\prime}-
by\os{y}\right)dx-2a\left[y\os{y}\right]_{0}^{1}\nonumber\\
&=&-2\int_{0}^{1}\os{y}\left(\DD{y}+2a\D{y}+by\right)dx\nonumber
\end{eqnarray}
The variational problem is thus equivalent to a Galerkin method applied to
the original problem with a weighting function which vanishes at the 
boundaries.
\section{A variational principle for Navier-Stokes equations}
We denote by $L^{2}(\Omega)$ the space of square integrable functions on 
$\Omega$ and by $H^{1}(\Omega)$ the Sobolev space made up of the functions 
which are in $L^{2}(\Omega)$. Let $H^{-1}(\Omega)=(H^{1}(\Omega))^{\prime}$ 
be the topological dual space of $H^{1}(\Omega)$ such that 
$H^{1}(\Omega)\subset L^{2}(\Omega)\subset H^{-1}(\Omega)$. The domain $\Omega$
is assumed to be simply-connected with a $C^{\infty}$ boundary 
$\partial\Omega$.

The Navier-Stokes equations for an incompressible flow, normalized to unit 
density, may be written in the form
\begin{eqnarray}
\left.\begin{array}{l}
\PDD{\B{u}}{t}-\nu\Delta\B{u}+\sum_{j=1}^{3}u_{j}\PDD{\B{u}}{x_{j}}+
\grad p=\B{f}\\
\dive\B{u}=0\\
\end{array}\right\}\hspace*{0.5cm}{\rm in}\,\,\Omega\cup\partial\Omega
\label{eq3}
\end{eqnarray}
with $\B{f}$ given in $H^{-1}(\Omega)$. 

We assume the boundary and initial conditions
\begin{eqnarray}
\begin{array}{ll}
\B{u}|_{\partial\Omega}=\B{U}, & \B{u}|_{t=0}=\B{a}\\
\end{array}\nonumber
\end{eqnarray}
where $\Delta$ is the Laplacian operator which is an isomorphism from 
$H^{1}(\Omega)$ onto $H^{-1}(\Omega)$ and $\nu$ is a positive parameter; the
kinematic viscosity.

We will now present a variational statement for the Navier-Stokes equations 
in the absence of any external forces $(\B{f}=0)$;
being a straightforward generalization of that for the one-dimensional damped
harmonic oscillator, with fixed end conditions, of the previous section. 

Consider the problem of finding stationary values for the functional
\begin{eqnarray}
J(u_{i},p,w_{i},r)=\int_{0}^{\tau}\Intt{\Omega}{\mathscr L}(u_{i},p,w_{i},r)
\,d\B{x}\,dt\label{eq3.1}
\end{eqnarray}
where ${\mathscr L}$ is the Lagrangian
\begin{eqnarray}
{\mathscr L}&=&\df{\nu}{2}\left(\PDD{u_{i}}{x_{j}}\right)^{2}+\df{1}{2}(w_{i}
+u_{i})u_{j}\PDD{w_{i}}{x_{j}}+u_{i}\PDD{p}{x_{i}}+\df{u_{i}}{2}
\PDD{w_{i}}{t}\nonumber\\
&-&\df{\nu}{2}\left(\PDD{w_{i}}{x_{j}}\right)^{2}-\df{1}{2}(u_{i}+w_{i})
w_{j}\PDD{u_{i}}{x_{j}}-w_{i}\PDD{r}{x_{i}}-\df{w_{i}}{2}
\PDD{u_{i}}{t}\nonumber
\end{eqnarray}
of the functions $u_{i}, p, w_{i}, r (i=1,2,3)$ dependent on the cartesian
spatial variables $\B{x}\in\Omega$ and the time interval $t\geq 0$. 
The double summation convention (summation over the reapeated indicies) is 
being used and the interval is over the domain $\Omega$ with smooth boundary 
$\partial\Omega$ and the time interval $[0,\tau]$. The class of admissible 
functions is such that 
\begin{description}
\item [(i)] $u_{i}$, $w_{i}$ have second order continuous spatial derivatives 
and first order continuous time derivative, $p$, $r$, have first order 
continuous spatial derivatives,
\item [(ii)] $u_{i}=w_{i}=f_{i}\,\,\forall\,\,\B{x}\in\partial\Omega$
for $t\in[0,\tau]$ where $\int_{\partial\Omega}f_{i}n_{i}=0$,
\item [(iii)] $u_{i}=w_{i}=g_{i}\,\,\forall\,\,\B{x}\in\Omega$ at
$t=0$ where $g_{i}$ are the components of a solenoidal vector field with 
which $f_{i}$ is compatiable. Hence $g_{i}$ is treated as an element of 
the space $H$ such that $g_{i}\in H^{\infty}(\Omega)\times
H^{\infty}(\Omega)\subset f_{i}(\Omega\rightarrow\Re^{3})$,
\item [(iv)] $u_{i}=w_{i}\,\,\forall\,\,\B{x}\in\Omega$ at $t=\tau$,
\item [(v)] $p=r\,\,\forall\,\,\B{x}\in\partial\Omega$ for $t\in[0,\tau]$.
\end{description} 
It will become clear that certain additional constraints are required when the
flow is steady and these will be stated at the appropriate time.

The variation of $J$ due to variations $\delta u_{i}$ etc. in the class of
admissible functions is
\begin{eqnarray}
\delta J(u_{i},p,w_{i},r)&=&\int_{0}^{\tau}\Intt{\Omega}\left[A_{i}(u,w,p)
\delta u_{i}-A_{i}(w,u,r)\delta w_{i}\right]\,d\B{x}\,dt\nonumber\\
&+&\int_{0}^{\tau}\Intt{\Omega}\left[B(u,w,p,\delta u,\delta w,\delta p)
-B(w,u,r,\delta w,\delta u,\delta r)\right]\,d\B{x}\,dt\label{eq3.2}
\end{eqnarray}	
where
\begin{eqnarray}
A_{i}(u,w,p)&=&\df{1}{2}\left[u_{j}\PDD{w_{i}}{x_{j}}+(w_{j}+u_{j})
\PDD{w_{j}}{x_{i}}\right]+\PDD{p}{x_{i}}+\df{1}{2}\PDD{w_{i}}{t}-\df{1}{2}
w_{j}\PDD{u_{i}}{x_{j}}\nonumber\\
B(u,w,p,\delta u,\delta w,\delta p)&=&\nu\PDD{u_{i}}{x_{j}}\PDD{\delta u_{i}}
{x_{j}}+\df{1}{2}(u_{i}+w_{i})u_{j}\PDD{\delta w_{i}}{x_{j}}+u_{i}
\PDD{\delta p}{x_{i}}+\df{u_{i}}{2}\PDD{\delta w_{i}}{t}\nonumber
\end{eqnarray}
The first integral involving $B$ may be integrated by parts to give
\begin{eqnarray}
& &\int_{0}^{\tau}\Intt{\Omega}B(u,w,p,\delta u,\delta w,\delta p)\,d\B{x}\,dt
=\nonumber\\
&-&\int_{0}^{\tau}\Intt{\Omega}\left\{\nu\PDDT{u_{i}}{x_{j}}\delta u_{i}+
\df{1}{2}\left[\PDD{(u_{i}+w_{i})u_{j}}{x_{j}}+\PDD{u_{i}}{t}\right]\delta 
w_{i}+\PDD{u_{i}}{x_{i}}\delta p\right\}\,d\B{x}\,dt\nonumber\\
&+&\int_{0}^{\tau}\Intt{\partial\Omega}\left\{\nu\PDD{u_{i}}{x_{j}}\delta u_{i}
+\df{1}{2}(u_{i}+w_{i})u_{j}\delta w_{i}+u_{j}\delta p\right\}n_{j}\,dS\,dt
+\Intt{\Omega}\left[\df{1}{2}u_{i}\delta w_{i}\right]_{0}^{\tau}\,d\B{x}
\label{eq3.3}
\end{eqnarray}

Conditions (ii) and (v) ensure that the boundary integral vanishes and
conditions (iii) and (iv) that the time independent integral vanishes. The
resulting Euler-Lagrange equations are
\begin{eqnarray}
\left.\begin{array}{ll}
p: & \PDD{u_{i}}{x_{i}}=0 \\
\\
r: & \PDD{w_{i}}{x_{i}}=0 \\
\\
u_{i}: & \nu\PDDT{u_{i}}{x_{j}}-\PDD{p}{x_{i}}=\PDD{w_{i}}{t}+
\df{u_{j}+w_{i}}{2}\left(\PDD{w_{i}}{x_{j}}+\PDD{w_{j}}{x_{i}}\right) \\
\\
w_{i}: & \nu\PDDT{w_{i}}{x_{j}}-\PDD{r}{x_{i}}=\PDD{u_{i}}{t}+
\df{w_{j}+u_{i}}{2}\left(\PDD{u_{i}}{x_{j}}+\PDD{u_{j}}{x_{i}}\right) \\
\end{array}\right\}\label{eq3.4}
\end{eqnarray}

At this point it should be noted that if $u_{i}=w_{i}$ then $\grad p$ and
$\grad r$ are identical and that there are only four independent equation;
these are identical in form to equations (\ref{eq3}) provided $p+u_{i}^{2}/2$
is interpreted as the pressure.

The difference functions
\begin{eqnarray}
\oss{v}{i}=\df{u_{i}-w_{i}}{2},\hspace*{1cm}\os{q}=\df{p-r}{2}\label{eq3.5}is not too large
\end{eqnarray}
satisfy
\begin{eqnarray}
\PDD{\oss{v}{i}}{x_{i}}&=&0\nonumber\\
\nu\PDDT{\oss{v}{i}}{x_{j}}-\PDD{\os{q}}{x_{i}}&=&-\PDD{\oss{v}{i}}{t}-
\df{u_{j}+w_{j}}{2}\left(\PDD{\oss{v}{i}}{x_{j}}+\PDD{\oss{v}{j}}{x_{i}}
\right)\label{eq3.6}
\end{eqnarray} 
subject to the boundary and initial conditions
\begin{eqnarray}
\begin{array}{lllll}
 & \oss{v}{i}=0 & \forall\B{x}\in\partial\Omega & \mbox{for} & t\in[0,\tau]\\
\mbox{and} & \oss{v}{i}=0 & \forall\B{x}\in\Omega & \mbox{at} & t=0,\tau\\ 
\end{array}\nonumber
\end{eqnarray} 
A problem essentially identical to that posed by equation (\ref{eq3.6}) arises
in providing uniqueness for weak solutions of the Navier-Stokes equations, see
for example, Ladyzhenskaya [3]. However, a more direct approach proceeds as
follows.

Multiplying (\ref{eq3.6}) by $\oss{v}{i}$ and using the incompressibility
conditions, we obtain
\begin{eqnarray}
\df{1}{2}\PDD{\oss{v}{i}^{2}}{t}=\nu\left(\PDD{\oss{v}{i}}{x_{j}}\right)^{2}+
\PD{x_{j}}\left\{\oss{v}{i}\os{q}-\nu\oss{v}{j}\PDD{\oss{v}{j}}{x_{i}}-
\df{u_{i}+w_{i}}{4}\oss{v}{j}^{2}\right\}-\df{u_{j}+w_{j}}{2}\oss{v}{i}
\PDD{\oss{v}{j}}{x_{i}}\label{eq3.7}
\end{eqnarray}
Integration of (\ref{eq3.7}) over the domain $\Omega$ and use of the boundary
condition $\oss{v}{i}=0$ on $\partial\Omega$ gives
\begin{eqnarray}
\df{1}{2}\df{dE}{dt}=\df{1}{2}\df{d}{dt}\Intt{\Omega}\oss{v}{i}^{2}\,d\B{x}=
\Intt{\Omega}\left[\nu\left(\PDD{\oss{v}{i}}{x_{j}}\right)^{2}-
\df{\oss{v}{i}\oss{v}{j}}{2}\PDD{(u_{j}+w_{j})}{x_{i}}\right]\,d\B{x}
\label{eq3.8}
\end{eqnarray}
As both $u_{i}$ and $w_{i}$ are components for solenoidal velocity fields the 
last term exeeds $-m\int_{\Omega}\oss{v}{i}^{2}\,d\B{x}$, where $m>0$ is the
maximum eigenvalue of the deformation tensor for the forcing flow,
\begin{eqnarray}
\df{1}{2}\left[\PDD{(u_{j}+w_{j})}{x_{i}}+\PDD{(u_{i}+w_{i})}{x_{j}}\right]
\nonumber
\end{eqnarray}
Integration of (\ref{eq3.8}) over the range $t_{0}<t<\tau$ where 
$0\leq t_{0}\leq\tau$ gives $E(t_{0})e^{2mt_{0}}<E(\tau)e^{2m\tau}$. Condition
(iii) has $u_{i}=w_{i}$ at $t=0\,\,\forall\,\,\B{x}\in\Omega$ thus $E(0)=0$ and
consequently $\oss{v}{i}=0$ and $\grad q=0\,\,\forall\,\,t\in[0,\tau]\,\,
\forall\,\,\B{x}\in\Omega$.

It is of some interest to discuss the precise nature of the stationary point 
for the functional. The following argument shows that the stationary point is 
indeed unique. Suppose, otherwise, that there are two stationary points 
corresponding to $(u_{i}^{1},p^{1},w_{i}^{1},r^{1})$ and
$(u_{i}^{2},p^{2},w_{i}^{2},r^{2})$. Each of these sets of functions satisfy 
the system (\ref{eq3.4}). The previous method may be applied to the differences
\begin{eqnarray}
u_{i}^{*}=\df{u_{i}^{1}-u_{i}^{2}}{2},\hspace*{1cm}p^{*}=\df{p^{1}-p^{2}}{2}
\nonumber
\end{eqnarray}
or
\begin{eqnarray}
w_{i}^{*}=\df{w_{i}^{1}-w_{i}^{2}}{2},\hspace*{1cm}r^{*}=\df{r^{1}-r^{2}}{2}
\nonumber
\end{eqnarray}
with the results that $u_{i}^{*}=w_{i}^{*}=0$ and $\grad p^{*}=\grad r^{*}=0$.
The functional is zero at the stationary point and this value is also taken at
any point where the test functions are such that $u_{i}=w_{i}$ and $p=r$. Now
the interchange of $(u,p)$ with $(w,r)$ changes the sign of $J$ and thus is 
clear that the stationary point is neither a maximum nor a minimum.
\section{A variational principle for the steady flow}
When the flow is time independent the appropriate functional is a suitably
modified version of (\ref{eq3.1}), i.e.,
\begin{eqnarray}
J^{s}(u_{i},p,w_{i},r)=\Intt{\Omega}{\mathscr L}(u_{i},p,w_{i},r)\,d\B{x}
\label{eq3.9}
\end{eqnarray}
Now the functions $u_{i}, p, w_{i}, r$ are required to be sufficiently smooth 
as in condition (i) and such that the pairs $(u_{i},p), (w_{i},r)$ take the
same values on the boundary $\partial\Omega$. A similar derivation for the
Euler-Lagrange equations leads to a time independent version of (\ref{eq3.4})
and the difference functions satisfy the steady versions of (\ref{eq3.6}) and 
the same boundary conditions. To show the uniquness of the stationary point in
non-steady problem, we set equation (\ref{eq3.8}) to zero. We can then see 
that 
\begin{eqnarray}
\nu\Intt{\Omega}\left(\PDD{\oss{v}{i}}{x_{j}}\right)^{2}\,d\B{x}&=&\left|
\Intt{\Omega}\df{\oss{v}{i}\oss{v}{j}}{2}\PDD{(u_{j}+w_{j})}{x_{i}}\,d\B{x}
\right|\nonumber\\
&\leq&\left\{\Intt{\Omega}\left(\oss{v}{i}\oss{v}{j}\right)^{2}d\B{x}\right\}
^{\tf{1}{2}}\left\{\Intt{\Omega}\df{1}{4}\left(\PDD{(u_{j}+w_{j})}{x_{i}}
\right)^{2}d\B{x}\right\}^{\tf{1}{2}}\nonumber
\end{eqnarray}
using the Schwarz's inequality, and in three-dimensions the above inequality 
is replaced by its generalized form given by Serrin [4]
\begin{eqnarray}
\nu\Intt{\Omega}\left(\PDD{\oss{v}{i}}{x_{j}}\right)^{2}d\B{x}\leq
3^{-3/4}\left\{\Intt{\Omega}\oss{v}{i}^{2}\,d\B{x}\right\}^{\tf{1}{4}}
\left\{\Intt{\Omega}\left(\PDD{\oss{v}{i}}{x_{j}}\right)^{2}d\B{x}\right\}
^{\tf{3}{4}}\left\{\df{1}{4}\Intt{\Omega}\left(\PDD{(u_{j}+w_{j})}{x_{i}}
\right)^{2}d\B{x}\right\}^{\tf{1}{2}}\nonumber
\end{eqnarray}
Payne and Weinberger [5] have provided a further inequality on the first 
integral on the right so that finally
\begin{eqnarray}
\nu\Intt{\Omega}\left(\PDD{\oss{v}{i}}{x_{j}}\right)^{2}d\B{x}\leq
\df{3^{-3/4}}{\lambda^{1/4}}
\Intt{\Omega}\left(\PDD{\oss{v}{i}}{x_{j}}\right)^{2}d\B{x}
\left\{\df{1}{4}\Intt{\Omega}\left(\PDD{(u_{j}+w_{j})}{x_{i}}
\right)^{2}d\B{x}\right\}^{\tf{1}{2}}\nonumber
\end{eqnarray}
where $\lambda$ is approximately $20R^{-2}$, and $R$ is the radius of the 
smallest sphere enclosing $\Omega$. The solution for $\oss{v}{i}$ is zero
when
\begin{eqnarray}
\df{R^{1/2}}{\nu}\left\{\df{1}{4}\Intt{\Omega}\left(\PDD{(u_{j}+w_{j})}
{x_{i}}\right)^{2}d\B{x}\right\}^{\tf{1}{2}}\leq R^{1/2}\lambda^{1/4}
3^{3/4}\approx 3\label{eq4.1}
\end{eqnarray}
Estimates for the left hand side of this expression in terms of boundary
data can be provided in the following way.

The steady state versions of equations (\ref{eq3.4}) imply that
\begin{eqnarray}
\df{\nu}{2}\PDDT{(u_{i}+w_{i})}{x_{j}}-\df{1}{2}\PDD{(p+r)}{x_{i}}
=\df{1}{4}(u_{j}+w_{j})
\left\{\PDD{(u_{i}+w_{i})}{x_{j}}+\PDD{(u_{j}+w_{j})}{x_{i}}\right\}
\label{eq4.2}
\end{eqnarray}
If equations (\ref{eq4.2}) are multiplied by $(u_{i}+w_{i})/2$ and then
integrated over the domain $\Omega$, we obtain
\begin{eqnarray}
& &\df{\nu}{4}\Intt{\Omega}\left[\PDD{(u_{i}+w_{i})}{x_{j}}\right]^{2}d\B{x}+
\df{1}{4}\Intt{\partial\Omega}\left[p+r+\df{(u_{j}+w_{j})^{2}}{4}\right]
(u_{i}+w_{i})n_{i}\,dS\nonumber\\
&=&\df{1}{4}\Intt{\Omega}(u_{i}+w_{i})(u_{j}+w_{j})
\PDD{(u_{j}+w_{j})}{x_{i}}\,d\B{x}\nonumber
\end{eqnarray}
The integral on the right hand side can be shown to be smaller than
\begin{eqnarray}
\df{3^{-3/4}}{\lambda^{1/4}}\left\{\df{1}{4}\Intt{\Omega}\left(
\PDD{(u_{j}+w_{j})}{x_{i}}\right)^{2}d\B{x}\right\}^{\tf{3}{2}}\nonumber
\end{eqnarray}
by a similar manipulation to that used previously so that if condition
(\ref{eq4.1}) is satisfied then
\begin{eqnarray}
& &\df{1}{4}\left\{\nu-\df{3^{-3/4}}{4\lambda^{1/4}}\Intt{\Omega}\left(
\PDD{(u_{j}+w_{j})}{x_{i}}\right)^{2}\,d\B{x}\right\}^{\tf{1}{2}}\Intt{\Omega}
\left(\PDD{(u_{i}+w_{i})}{x_{j}}\right)^{2}d\B{x}\nonumber\\
&<&\left|\df{1}{4}
\Intt{\partial\Omega}\left[p+r+\df{(u_{j}+w_{j})^{2}}{4}\right](u_{i}+w_{i})
n_{i}\,dS\right|\nonumber
\end{eqnarray}
As in the non-steady problem, the stationary point is unique and is neither
a maximum or a minimum.
\section{An extended variational principle}
It is possible to remove a number of the restrictions on the class of 
admissible functions by adding a surface integral to the previous functional
$J$ given by (\ref{eq3.1}). Consider the new functional
\begin{eqnarray}
I=J+\int_{0}^{\tau}\Intt{\partial\Omega}{\mathscr M}_{j}(u_{i},p,w_{i},r)n_{j}\,
dS\,dt\label{eq7.1}
\end{eqnarray}
where $n_{j}$ are the components of the outward normal to $\partial\Omega$ and
\begin{eqnarray}
{\mathscr M}_{j}(u_{i},p,w_{i},r)&=&-u_{j}^{s}p+w_{j}^{s}r+\df{u_{j}u_{i}^{2}+
u_{i}^{2}w_{j}}{4}-\df{w_{j}w_{i}^{2}+w_{i}^{2}u_{j}}{4}\nonumber\\
&+&\nu\left[-(u_{i}-u_{i}^{s})\PDD{u_{i}}{x_{j}}+(w_{i}-w_{i}^{s})\PDD{w_{i}}
{x_{j}}\right.\nonumber\\
&+&\left.\eps_{ijk}(w_{k}-w_{k}^{s}-u_{k}+u_{k}^{s})\PD{x_{i}}|\B{u}-\B{u}^{s}
+\B{w}-\B{w}^{s}|\right]\nonumber
\end{eqnarray}
where $\eps_{ijk}$ is the permutation tensor and the superscript $s$ indicates
the surface values.

The derivation of this expression for ${\mathscr M}_{j}$ is complicated by the 
fact that when $\nu=0$ the physical boundary condition on $\partial\Omega$
should relate only to the normal velocity component: accordingly, normal and
tangential components require separate treatment.

We require $I$ to be stationary $(\delta I=0)$ subject to the conditions (i),
(ii) and (iv) where in addition $u_{i}=w_{i}$ at one point on $\partial\Omega$
and the variation are such that $\delta u_{i}+\delta w_{i}=0$ on
$\partial\Omega$. It is then readily verified that variations of $u_{i}, p,
w_{i}, r$ leads to the Navier-Stokes equations together with its adjoint
system.

However, the boundary conditions on $\partial\Omega$ require a little care. 
Incorporating the additional integral of ${\mathscr M}_{j}$, the variations in
$w_{j}$ gives
\begin{eqnarray}
n_{i}(u_{i}-u_{i}^{s})(u_{j}-u_{j}^{s})+\nu n_{i}\eps_{ijk}\PD{x_{k}}|\B{u}-
\B{u}^{s}|=0\label{eq6.1}
\end{eqnarray}
Taking the scalar product of (\ref{eq6.1}) with $n_{j}$, we obtain
\begin{eqnarray}
n_{i}(u_{i}-u_{i}^{s})=0\,\,\,\,\forall\,\,\,\,\B{x}\in\partial\Omega\nonumber
\end{eqnarray}
independently of whether $\nu$ is zero or not. This is, of course, the 
inviscid boundary conditions, which prescribes the normal velocity component
at each point on $\partial\Omega$. Accordingly, the second term in 
(\ref{eq6.1}) must itself vanish. In vector notation, this is just 
$\nu(\B{n\times}\grad|\B{u}-\B{u}^{s}|)=0$, which if $\nu\neq 0$, asserts that
$|\B{u}-\B{u}^{s}|$ remains constant on $\partial\Omega$. If we now introduce
the further requirement that $\B{u}=\B{u}^{s}$ at a {\em single point} on
$\partial\Omega$, we have the boundary conditions for viscous flow that
\begin{eqnarray}
u_{i}=u_{i}^{s}\,\,\,\,\forall\,\,\,\,\B{x}\in\partial\Omega\nonumber
\end{eqnarray}
Variations of $p$ and $u_{i}$ give the corresponding boundary conditions for
$w_{i}$ on $\partial\Omega$. Variation of $p$ leads immediately to
\begin{eqnarray}
n_{i}w_{i}=0\label{eq6.2}
\end{eqnarray}
corresponding to the inviscid boundary condition for $u_{i}$. Variation of 
$u_{i}$, subject to the requirement that $\delta u_{i}=0$ on $\partial\Omega$,
leads to
\begin{eqnarray}
\nu n_{i}w_{j}\left\{\PDD{\delta u_{i}}{x_{j}}+\PDD{\delta u_{j}}{x_{i}}-
\eps_{ijk}\PDD{|\delta\B{u}|}{x_{k}}\right\}=0\,\,\,\,\forall\,\,\,\,\B{x}\in
\partial\Omega\nonumber
\end{eqnarray}
However, since $\delta u_{i}=0$ on $\partial\Omega$,
\begin{eqnarray}
\PDD{\delta u_{i}}{x_{j}}&=&n_{i}\PDD{\delta u_{i}}{\ell}\nonumber\\
\PDD{|\delta\B{u}|}{x_{k}}&=&n_{k}\PDD{|\delta\B{u}|}{\ell}\nonumber
\end{eqnarray}
where $\ell$ denotes distance along the outward normal on $\partial\Omega$.
This yields
\begin{eqnarray}
\nu(n_{i}n_{j}w_{j}+w_{i})\PDD{\delta u_{i}}{\ell}=0\nonumber
\end{eqnarray}
the third term vanishing identically. Using (\ref{eq6.2}) we have, from 
variation in $u_{i}$,
\begin{eqnarray}
w_{i}=0\,\,\,\,\forall\,\,\,\,\B{x}\in\partial\Omega\nonumber
\end{eqnarray}
whenever $\nu\neq 0$. There is of course a corresponding result for a steady
state problem.
\section{Summary}
To summarise the results obtained in this work we note the following:
\begin{description}
\item [(a)] the functional (\ref{eq3.1}) of the eight functions $u_{i}, p,
w_{i}, r (i=1,2,3)$ which belong to the class of functions (i)-(v) has a unique 
stationary value of zero at the point $(v_{i}, q, v_{i}, q)$ where the 
equations satisfied by $v_{i}$ and $q$ are the Navier-Stokes equations
\begin{eqnarray}
\PDD{v_{i}}{x_{i}}&=&0\nonumber\\
\PDD{v_{i}}{t}+v_{j}\PDD{v_{i}}{x_{j}}&=&-\PD{x_{i}}\left(q+\tf{1}{2}v_{j}^{2}
\right)+\nu\PDDT{v_{i}}{x_{j}}\nonumber
\end{eqnarray}
within a domain $\Omega$ with smooth boundary $\partial\Omega, v_{i}$ is given
on $\Omega$ at $t=0$ and on $\partial\Omega$ for $t\in[0,\tau]$.
\item [(b)] The functional (\ref{eq3.9}) of the eight functions $u_{i}, p, 
w_{i}, r (i=1,2,3)$ which belong to the class functions satisfying
\begin{description}
\item [(i)] $u_{i}, w_{i}$ have continuous second order spatial derivatives
and $p$ and $r$ have continuous first order spatial derivatives within the
domain $\Omega$,
\item [(ii)] $u_{i}=w_{i}=g_{i}\,\,\forall\,\,\B{x}\in\partial\Omega$ where
$\int_{\partial\Omega}g_{i}n_{i}\,dS=0$,
\item [(iii)] $p=r\,\,\forall\,\,\B{x}\in\partial\Omega$,
and\item [(iv)] condition (\ref{eq4.1}) is satisfied,
\end{description}
has a unique stationary value zero at the stationary point $(v_{i}, q, v_{i},
q)$ where $v_{i}, q$ satisfy the steady state Navier-Stokes equations
\begin{eqnarray}
\PDD{v_{i}}{x_{i}}&=&0\nonumber\\
v_{j}\PDD{v_{i}}{x_{j}}&=&-\PD{x_{i}}\left(q+\tf{1}{2}v_{j}^{2}\right)+
\nu\PDDT{v_{i}}{x_{j}}\nonumber
\end{eqnarray}
within a domain $\Omega$ with smooth boundary $\partial\Omega$ and where 
$v_{i}$ is given on $\partial\Omega$.
\end{description}
It is to be noted that in the steady state case if the Reynolds number does 
not remain finite (\ref{eq4.1}) will not be satisfied and so its use is
somewhat restricted.  

Although these variational integrals yield the correct equations in $\Omega$,
they do not incorporate the correct physical boundary conditions on
$\partial\Omega$. To circumvent this difficulty, it is necessary to add a
further surface integral, namely (\ref{eq7.1}).

\section*{References}
\begin{description}
\item {[1]} Millikan, C.B. (1929) On the steady motion of viscous 
incompressible fluids, with particular reference to a variational principle.
{\em Phil. Mag.}, {\bf 7}, 641.
\item {[2]} Finlayson, B.A. (1972) Existence of variational principles for
the Navier-Stokes equations. {\em Phys. Fluids}, {\bf 15}, 963.
\item {[3]} Ladyzhenskaya, O.A. (1969) {\em The mathematical theory of viscous
incompressible flow}. Gordon and Breach.
\item {[4]} Serrin, J. (1963) The initial value problem for the Navier-Stokes
equations. In {\em Nonlinear problems} ed. R.E. Langer, The University of 
Wisconsin press, page 98.
\item {[5]} Payne, L.E and Weinberger, H.F. (1963) An exact stability bound for
Navier-Stokes flow in a sphere. In {\em Nonlinear problems} ed. R.E. Langer, 
The University of Wisconsin press, page 311.
\end{description}

\end{document}